\documentclass[10pt]{article}

\usepackage{graphicx}
\usepackage{graphicx,psfrag}
\usepackage{amsmath,amssymb,amsthm}
 
\def \G {\Gamma}

\def \vs {\vskip 0.2cm}

\def \n {\noindent}

\def \E {{\bf  \, E  }}

\def \ch  {{\hbox{\hskip 0.03cm ch}}}
\def \sh {{\hbox{\hskip 0.03cm sh}}}

\def \Tr {{\mathrm{\, Tr}}}

\def \bC {{\mathbb C}}

\def \bE {{\mathbb E}}

\def \bN {{\mathbb N}}

\def \bR {{\mathbb R}}

\def \CB{{\cal B}}

\def \CE {{\cal E}}

\def \CN {{\cal N}}
\def \CT {{\cal T}}

\def \CI {{\cal I}}

\def \CP {{\cal P}}

\def \CV {{\cal V}}

\def \G{{\Gamma}}

\def \a {\alpha}

\def \l {\lambda}

\def \s {\sigma}
\def \t  {\tau }

\def \th {\theta}

\def \vep {\varepsilon}

\begin{document}
\title{
On asymptotic properties of high moments of compound Poisson 
distribution
}

\author{O. Khorunzhiy\\ Universit\'e de Versailles - Saint-Quentin \\45, Avenue des Etats-Unis, 78035 Versailles, FRANCE\\
{\it e-mail:} oleksiy.khorunzhiy@uvsq.fr}
\maketitle
\begin{abstract}
We study    moments $M_k(\lambda)$ of the sum of 
 random variables $X_1+\dots +X_{N_\l}$,
where $N_\l$ follows the Poisson probability distribution   with  mean value 
$\l$ and $\{X_j\}$ is a family of i.i.d. random variables also independent from $N_\l$. We obtain 
an explicit expression for the leading term of the asymptotic expansion of $M_k(x)$ as $k\to\infty$. 
  We show that if $\l_k$ is much smaller or proportional to $k$, then 
  asymptotic behavior of $M_k(\l_k)$ 
  is determined by the exponential generating function of 
     $X_j$
   while 
  in the asymptotic regime when $\l_k$ is much greater than  $k$, the leading term of 
$k^{-1}\ln M_k(\l_k)- \ln \l_k$, $k\to\infty$ depends on the first non-zero moment of $X_j$ only.

As a consequence, 
we establish a concentration property 
of  maximal vertex degree of large weighted random graphs. 
Another  application is related with a variable that arises in the studies
of high moments of large random matrices. 
Finally, regarding three  particular  cases of probability distribution of 
$X_j$, we comment on the  asymptotic  behavior of  
certain combinatorial polynomials, including  the Bell 
polynomials of even partitions. 

  \end{abstract}



\section{Introduction, main result and discussion}
Compound Poisson distribution is widely used
in a number of applications in \mbox{various} areas, with the majority of recent applications in financial and risk modeling  
(see monograph \cite{Ross} and references therein). 
 This distribution can be \mbox{associated} with  a random variable 
$$
Y_\l   =\sum_{j=1}^{N_\l } X_j,
\eqno (1.1)
$$
where $\{X_j\}_{j\in \bN}$ is a family of i.i.d.  random variables and random variable $N_\l$  independent from $\{X_j\}$ 
follows the Poisson law  with mean value $\l$
\cite{KPST}. 
We denote by $\bE$ the mathematical expectation with respect to a measure generated 
by the family of random variables $\left\{ N_\l, \{X_j\}_{j\in \bN}\right\}$.  
We assume that all moments of $X_j$ exist and denote $V_i= \bE X_j^i$, 
$i\in \bN$.

It is known  that the moments of $Y_\l  $ are given by the following expression, 
$$
M_k(\l )= \bE Y_\l ^k = k! \sum_{(l_1, l_2, \dots, l_k)^*}^k 
\prod_{i=1}^k {(\l  V_i)^{l_i}\over (i!)^{l_i} l_i!}, \quad k\in \bN, 
\eqno (1.2)
$$
where the sum runs  over  all  $l_i\ge 0$ such that 
$l_1+2l_2 + \dots + kl_k=k$.
In the case when the random variables $X_j$ are all equal to $1$, relation (1.2) 
determines  the one-variable Bell polynomials $\CB_k(x)$ \cite{Carl},
$$
M_k^{(V_i=1)}(x)= \CB_k(x) = k! \sum_{(l_1, l_2, \dots, l_k)^*}^k 
\prod_{i=1}^k {x^{l_i}\over (i!)^{l_i} l_i!}, \quad x\in \bR.
 \eqno (1.3)
 $$ 
If $x=1$, then the right-hand side of (1.3) coincides with  the $k$-th Bell number $\CB_k(1) = B^{(k)}$  \cite{Bell-1,Bell-2,Berndt,FS,Touch},
$$
B^{(k)} = \sum_{(l_1, l_2, \dots, l_k)^*}^k B^{(k)}_{l_1, l_2, \dots, l_k}\, ,
\quad 
B^{(k)}_{l_1, l_2, \dots, l_k} = {k!\over (1!)^{l_1}l_1! \, (2!)^{l_2} l_2!
\cdots (k!)^{l_k} l_k!}.
\eqno (1.4)
$$
The Bell number $B^{(k)}$ counts the number of all possible partitions $\pi(k)$ 
of the set of $k$ elements into non-empty subsets (blocks). Then $l_i$ represents the number of blocks of size $i$ and the sum $l_1+\dots+ l_k= \vert \pi(k)\vert $
is the total number of blocks  in the partition $\pi(k)$.   
According to this, we can say that polynomials $M_k(x)$ (1.2)
represent a generalization of $\CB_k(x)$. This generalization can be referred as to    the  Bell polynomials of weighted partitions or simply as to the weighted Bell polynomials. 
\vs 
Limiting behavior  of the Bell numbers $B^{(k)}$ has been studied 
since 50-s \cite{FS,MoWy,TE}  (see also monograph \cite{FS})
while  questions related with  asymptotic properties of Bell polynomials
represent more recent  area of researches \cite{Do,Elb,Zhao}. 
In paper \cite{CR},  moments (1.2) of compound Poisson distribution have been studied 
and recurrent relations have been established for $M_k(\l )$ with given $k$.

   It should be noted that probabilities of large deviations for sums of random number of i.i.d. random variables  
 of the form (1.1) have been studied in the limit $\lambda\to\infty$ \cite{Mita}. However, 
 these results were obtained  with no reference to the moments 
 $M_k(\l)$ of the corresponding compound Poisson process. 
 The use the high moments $M_k(\l)$ with increasing $k$ 
 can be helpful also in studies  
 of the concentration properties of  sums $Y^{(i)}_\l$, $\l\to\infty$  (1.1) and their maxima that is in close relation  with the  maximal vertex degree of weighted random graphs.
Also,  the moments of  $Y_\l$ mimic certain elements arising in studies of  high moments of large random matrices of the Wigner ensemble.
Here, the limiting transition when  $k$ and $\l$ tend simultaneously to infinity  
is naturally motivated.
Recently, high moments of the compound Poisson probability distribution 
have been used as  upper bounds of cumulants of number of weighted walks
on large random graphs \cite{K-24}. 
 We briefly discuss some of these questions at the end of the paper.

Up to our knowledge,  asymptotic properties of  weighted versions  of Bell polynomials $M_k(x_k)$ with infinitely increasing $k$ have not been yet 
  considered and the main statement proved of this paper is new. 
   formulate our results under assumption that 
there exists finite or infinite $u_0>0$ such that
 the following 
series converges,
$$
\forall u \in [0, u_0): \ H(u)= \sum_{k=0}^\infty {V_{k}\over k!}  { u^{k}}<+\infty.
 \eqno (1.5)
$$
Function $H(u)$  is  known as the moment generating function of the probability distribution $P_X$, $H(u)= \bE e^{uX}$. 
  Clearly, all the derivatives  of $H(u)$  exist for all $u\in [0, u_0)$ and   the first derivative 
$H'(u)$ is strictly increasing on $(0,u_0)$ \cite{Baha}. 
In one of the three asymptotic regimes considered in this paper, we impose the following additional conditions on 
$H(u)$:
\vs 

i)  there exists $0< \alpha< 1/2$ such that 
$$
\hbox{either} \quad {H(u+1)\over (H''(u))^{1+\a}} \to 0, \ u\to \infty\quad 
{\hbox{or}} \quad 
{H((u_0+u)/2)\over (u_0 - u)^3  (H''(u))^{1+\a}} \to 0, \ u\to u_0
\eqno (1.6)
$$
and
$$
u^4 \left( H''(u)\right)^{3+2\a} \ll H(u)^4, \ \hbox{when either } u\to \infty \ \hbox{or} \ u\to u_0;
\eqno(1.7)
$$

\vs

ii) there exists  positive $w_0$ and $\beta$ such that 
$$
\bE e^{u W}{\hbox{I}}_{\{W\ge w_0\}} \ge \beta H(u).
\eqno (1.8)
$$
In our studies, we distinguish three major asymptotic regimes
  in dependence whether $x_k$ is much less, proportional or much greater than $k$ when $k$ tends to infinity.


\vs \n  {\bf Theorem 1.1.} 
\n {\it
Consider a sequence $(x_k)_{k\in \bN}$  such that $x_k = \chi k$,
$\chi>0$.
Let $u>0$ 
  be determined by equation
$$
u H'(u)= {1\over \chi}.
\eqno(1.9)
$$
Then the following asymptotic equality is true,
$$
M_k(x_k)= \left( { k\over u} \exp\left\{ {H(u)-1\over u H'(u)}-1\right\} \right)^k   {1+o(1)\over \sqrt{\chi   u(H'(u) + uH''(u))}}
\, 
,
\quad k\to \infty.
\eqno (1.10)
$$

If $(x_k)_{k\in \bN}$ is such that $x_k = \chi_k k$ with $\chi_k\to\infty$ as $k\to\infty$,
then relation (1.10) is true with $u$ replaced by $u_k$ determined by 
equation 
$$
u_k H'(u_k)= {1\over \chi_k}, \quad k\in \bN.
\eqno(1.11)
$$

If   conditions (1.6), (1.7) and (1.8)  are satisfied, then
relation (1.10) remains valid in the limiting transition 
 $k\to\infty$ such that $\chi_k\to 0$
with $x_k$ finite or tending to infinity  and  $(u_k)_{k\in \bN}$ determined by (1.11). }

\vs 
{\it Remark.} Relation (1.10) remains valid also in the limit with vanishing $x_k$
 such that the following version of (1.6) holds,
$$
\hbox{either} \quad {H(u_k+1)\over \sqrt {x_k} (H''(u_k))^{1+\a}} \to 0
\quad 
{\hbox{or}} \quad 
{H((u_0+u_k)/2)\over\sqrt {x_k} (u_0 - u_k)^3  (H''(u_k))^{1+\a}} \to 0, \quad k\to\infty
$$
with $u_k$  still determined by (1.11).

\vs
Theorem 1.1 generalizes  results obtained  in  \cite{K-1} for the one-variable Bell polynomials (1.3) with constant weight $X_j=1$. 
To better see  the novelties   due to $X_j$, we present  the following    three corollaries of Theorem 1.1.

\vs 
\n {\bf Corollary  A.} Let us consider  first the asymptotic regime when $1\ll k\ll x_k$. In this 
case
$\chi_k \to\infty$ and we denote this limiting transition by
$$
(k,x_k)_\infty \to \infty.
\eqno (1.12)
$$

It is not hard to deduce from (1.9) that in this case  $u_k\to 0$ and Theorem 1.1 implies that
$$
M_k(x_k) = \begin{cases}
 \left( x_k V_1 \big(1+o(1)\big)\right)^k ,  & \text{if  $V_1\neq 0$} , \\
\left( {\displaystyle x_k k V_2\over \displaystyle  e}\big(1+o(1)\big)\right)^{k/2}, \quad   & 
\text {if  $V_1=0$.} 
\end{cases}
\eqno (1.13)
$$
Let us note that if the probability distribution of $X_j$ 
is symmetric, then Theorem 1.1 holds for the even moments $M_{2k}(x)$ only.

 In the case of constant $X_j=1$, the first relation of (1.13) has been obtained in \cite{K-1} for Bell polynomials 
 (1.3). The second is obtained for the case of restricted Bell polynomials \cite{K-1} that concern the partitions (1.4) that have no blocks of size one.  
\vs 
\vs

\n {\bf Corollary  B.}  Let us consider the second asymptotic regime  when   
$x_k/k \to  \chi  $ as $k\to\infty$. We denote this limiting transition by
$$
(k,x_k)_\chi\to\infty.
\eqno (1.14)
$$
It follows from Theorem 1.1 that  in this asymptotic regime
$$
M_k(x_k) = \left(x_k e^{\Psi(\chi)}\big(1+o(1)\big)\right)^k, \quad 
(k,x_k)_\chi \to \infty,
\eqno (1.15)
$$
where 
$$
\Psi(\chi) =
  {H(u)-1\over u H'(u)}- 1  + \ln H'(u),
\eqno (1.16)
$$
and  $u$ is given by a solution of equation (1.9).

\vs This result can be also reformulated as follows,
$$
M_k(x_k)= \left( {k\over e} A_\chi \big( 1+o(1)\big)\right)^k, \quad (x,u)_\chi\to
\infty,
\eqno (1.17)
$$
where 
$$
A_\chi = {1\over u} \exp\left\{ {H(u)-1\over uH'(u)}\right\}
\eqno (1.18)
$$ 
and $u$ is determined by (1.9). Relations (1.13) can be deduced from (1.15) and (1.16) in the limit $u\to 0$.

\vs Let us note  that in the case of $X_j=1$
we have  obviously $H(u)=e^u$
and 
relation (1.9) takes the form of  the Lambert equation \cite{dB}
$$
ue^u= t, \quad t=  \chi^{-1}.
\eqno (1.19)
$$
Then (1.15) is transformed into asymptotic equality 
$$
M_k^{(X_j=1)}(x_k)= \CB_k(x_k)= \left( x_k e^{\psi(\chi)} \big( 1+ o(1)\big)\right)^k,
\eqno (1.20)
$$
where 
$$
\psi(\chi)= {e^u-1\over u e^u} + u-1 .
\eqno (1.21)
$$
This result  has been previously obtained   with the help of the ray method applied to the differential-difference equation satisfied by $\CB_k(x)$
\cite{Do,Elb}. 
\vs

\noindent {\bf Corollary  C.}  Let us consider the third asymptotic regime when 
$x= o(k)$, $k\to\infty$. 
denote this limiting transition as
$$
(k,x_k)_0\to \infty.
\eqno (1.22)
$$
In this  regime,  
dependence of the asymptotic behavior of $M_k(x_k)$ on the probability distribution 
of $X_j$ manifests itself in the most pronounced way. 
In Section 3 we consider several key  particular cases of the probability distribution of $X$. 
Here we present a part of the results obtained.

In the case when $X_j$ have exponential 
 gamma distribution $\G{(m,\theta)}$ with $m>8$, the
asymptotic behavior  of the moments $M_k(x_k)$ is similar  to 
that given by  (1.18),
$$
M_k(x_k) = \left( {k\th \over e}\big(1+o(1)\big)\right)^k,  \quad  (k,x_k)_0\to\infty.
\eqno (1.23)
$$
In the case when $X_j$ follow  the  Gaussian (normal) law $\CN(a,\s^2)$,
we get asymptotic relation 
$$
M_k(x_k) = \left( { k \sigma \over e \sqrt {2  \big( \ln k - \ln x_k\big)}}
\big( 1+o(1)\big)\right)^k, \quad (k,x_k)_0\to \infty.
\eqno (1.24)
$$
In the case of the standard centered normal distribution, $X_j \sim N(0,1)$, even moments show the following asymptotic behavior,
$$
M_{2k}(x_k) = \left( {\sqrt 2 k\over e \, \sqrt{ \ln (2k) -\ln x_k}}
\big(1+o(1)\big)\right)^{2k}, \quad (k,x_k)_0\to \infty.
\eqno (1.25)
$$

The wighted Bell polynomials with weights given by centered Bernoulli random variables show the following asymptotic behavior,
$$
M_{2k}(x_k) = \left( {2k\over e \big(\ln (2k) - \ln x_k\big) }
\big(1+o(1)\big)\right)^{2k}, \quad (k,x_k)_0\to \infty.
\eqno (1.26)
$$
This asymptotic expression coincides with that obtained for 
 the pure (non-weighted) Bell polynomials,
$$
\CB_k(x_k) = \left( {k\over e\big(\ln k - \ln x_k \big)}(1+o(1))\right)^k,
\quad (k,x_k)_0\to \infty
\eqno (1.27) 
$$
with $k$ being replaced by $2k$ \cite{Do,Elb,K-1}.

\vs

Regarding  corollaries A, B and C of Theorem 1.1, we see that 
 in the last asymptotic regime (1.22) when $x_k\ll k$, the properties of the weights $X_j$ crucially modify the asymptotic behavior of the moments $M_k(x_k)$. 
In the second asymptotic regime (1.14) when $x_k$ is proportional to $l$, the moments $M_k(x_k)$ exhibit almost the same asymptotic behavior (1.18), where the form of the probability distribution of $X_j$ enters into  the right-hand side of (1.18) via a constant $A_\chi$. 
The same is true for the expression (1.15). Finally, we see that in the asymptotic regime (1.12) when $x_k\gg k$, the moments $M_k(x_k)$ are not sensible to the details of  the  probability distribution of the weights $X_j$ and the leading term of their
asymptotic expansion is almost universal. 
This universality can be explained by a kind of the law of large numbers
that follow random variables $Y_\l/\l$  
in the limit $\l \to\infty$. We discuss this topic in Section 4.

 The paper is organized as follows. In Section 2, we formulate and prove our main technical result  
given by  the Local Limit Theorem for a family of  auxiliary random variables. 
Then we prove Theorem 1.1 as well as the corollaries A and B related with the  asymptotic regimes (1.12) and (1.14), respectively. 
In Section 3, we consider the asymptotic behavior of the moments $M_k(x_k)$ in  the third limiting transition (1.22). In Section 4, 
 we formulate and prove  auxiliary statements and discuss  a number of supplementary facts related to our main results.

\section{Proof of Theorem 1.1}

We prove Theorem 1.1 with the help of the method proposed in \cite{TE}
to study the asymptotic behavior of the Bell numbers $B^{(k)}$   and then modified  in \cite{K-1} in applications to the case of Bell polynomials $B_k(x)$ (1.3). This method is based on the observation
that the local limit theorem is valid for an auxiliary random variable 
$Z$ such that the probability $P(Z=k)$ is proportional to $ B_k$.
 The use of the local limit theorem in order to get the asymptotic properties of combinatorial items dates back to the works 
by E. A. Bender (see \cite{Be} and also \cite{FS,GLN}). 
Further use of this technique developed in \cite{Be} and other papers  would require 
proofs of more statements such as the log-concavity of 
the sequence $M_k(x)$. We stay within the frameworks
of the stochastic version of the method outlined in \cite{TE}.

\subsection{Random variables and Central Limit Theorem}

Let us introduce a random variable 
$Z^{(x,u)}$ that takes values in $\bN$ such that 
$$
P(Z^{(x,u)} = k)= M_k(x) {u^k\over k!\,  G(x,u)}, \quad j\in \bN, \ u>0,
\eqno (2.1)
$$
where 
$
G(x,u)
$
is the normalization factor. 
It is well known that this exponential generating function $G(x,u)$
is determined by relations
$$
G(x,u)= \sum_{j=0}^\infty M_j(x) {u^j\over j!}
 = 
\exp\left\{ x \big(H(u) - 1\big)\right\}.
\eqno (2.2)
$$
The last equality of (2.2) can be proved with the help of  standard combinatorial arguments. 
This equality relates the moments of $X$ and $M_k(x)$. 
The main task of this paper is to see in what way and in what extension 
$H(u)$ determines the asymptotic behavior of $M_k(x)$ in one or another asymptotic regime. 
It is not surprising that relation (2.2) represents a keystone of the method in general and of the proof of
our main results. 

It is easy to see that the generating function of the probability distribution (2.1) determined by  
$F_{x,u}(\t ) = \sum_{j=0}^\infty p_j^{(x,u)} \t ^j$ 
allows the following representation,
$$
F_{x,u}(\t ) = {G(x,\t u)\over G(x,u)}.
\eqno (2.3)
$$ 
Elementary calculations based on (2.1), (2.2) and (2.3) show that 
$$
\E Z^{(x,u)} = \sum_{j=0}^\infty j p_j^{(x,u)} = 
xu H'(u) 
\eqno (2.4)
$$
and 
$$
Var(Z^{(x,u)})= \s_{Z}^2 = x\big( uH'(u) + u^2 H''(u)\big).
 \eqno (2.5)
$$
\vs

The main technical statement  given by the Local Central Theorem presented below says  that if  $k$ is large and
not far from the mean value of $Z^{(x,u)}$ (2.4),  then 
the probability
$P(Z^{(x,u)}=k)$
is close to $({2\pi} \sigma_{Z}^2)^{-1/2}$,
$$
P(Z^{(x,u)}=k)=  { {1\over \sqrt{2\pi} \s_{Z}}}(1+o(1)), \quad \E Z^{(x,u)}=
k(1+o(1)), \quad  k\to \infty.
\eqno (2.6)
$$
Rewriting (2.1) in the form
$$
M_k(x)=  {k!\,\over u^k} \,  \exp\left\{x\big( H(u) - 1\big)\right\}  \, 
 P(Z^{(x,u)}=k)
 \eqno(2.7)
$$
and using  (2.6), on can obtain the  main result of Theorem 1.1 given by 
(1.10).

\vs 
It is interesting to note that 
 while the first part of  relation (2.6) represents the principal result
 of this  paper, the last factor of the right-hand side of (2.7), in the majority of cases, 
 does not play a crucial role  with respect to  the asymptotic behavior of $M_k(x)$. In contrast with this, 
the second part of (2.6) that determines the value of $u$ in dependence of $k$ and 
$x$ contributes essentially to the form of   asymptotic expansions of $M_k(x)$. 
We discuss these questions in more details in Section 3. 

Let 
$$
Y^{(x,u)}=  {Z^{(x,u)}-\E Z^{(x,u)} \over \s_Z}. 
\eqno (2.8) 
$$ 
We consider an infinite sequence of $(x_k, \chi_k)$ and determine $u_k$ by relation (1.9). We introduce a function
 $$
 \Phi_{Y_k}(t) = \E \exp\{ - it Y_k\},
 $$
 where $Y_k= Y^{(x_k,u_k)}$ and formulate our first  statement  given by  the Central Limit Theorem for the sequence of random variables $Y_k $, $k\to\infty$. 

\vs 
{\bf Lemma 2.1.} 
{\it For any given $t\in \bR$, the following asymptotic equality 
$$
 \Phi_{Y_k}(t)=   \exp\{ - t^2/2\}\big(1+o(1)\big)
\eqno (2.9)
 $$
 holds in the limiting transitions $(k,x_k)_\infty\to \infty$ (1.12) and   $ (k,x_k)_\chi\to\infty$ (1.14) ; 
 if $H(u)$ is such that  conditions (1.6) and  (1.8) are verified, then  
 relation (2.9) holds also in the limiting transition $(x,u)_0\to\infty$ (1.22).}

\vs 
{\it Proof of Lemma 2.1.}
Relations (2.1) and (2.2) obviously imply  that
$$
\Phi_Y(t) = \exp\{  - it \E Z/\s_Z\} F\big(e^{it/\s_Z}\big),
\eqno (2.10)
$$
where 
$$
F\big(e^{ it /\s_Z}\big)=  { G\big(x,u e^{it/\s_Z}\big)\over G(x,u)} = \exp\left\{ x\left(H
\big(ue^{it/\s_Z}\big) - H(u)\right)\right\}.
\eqno(2.11)
$$
Here and below we omit the subscript $k$ everywhere when no confusion can arise.

Taking into account (1.5), we conclude that the series 
$
H(z) = \sum_{k\ge 0} {V_k} \,  z^k/k! 
$ 
converges uniformly in the open disk $B(0, u_0) = \{ z\in \bC: \, \vert z \vert < u_0\}$ and therefore the following analog of the Taylor expansion holds,
$$
H(z) = \sum_{j=0}^2 {H^{(j)}(u)\over j!} \, (z-u)^j + R_2(z,u), \quad u\in [0, u_0)
\eqno (2.12)
$$
with 
$$
\vert R_2(z,u)\vert \le h(u, r) 
{ \big(\vert z-u\vert/r\big)^3\over 1 - \vert z- u\vert/r},
\quad h(u,r) = \max_{s\in \bC: \, \vert s- u\vert = r} \vert H(s)\vert,
\eqno (2.13)
$$
where one can take $r= (u_0-u)/2$ if $u_0<\infty$  and $r= 1$ if $u_0=+\infty$.
Relations (2.12) and (2.13) are proved in Section 4.

Regarding the right-hand side of (2.12) with $z = u e^{it/\s_Z}$
we conclude that
$$
H(u e^{it/\s_Z}) - H(u) = uH'(u)\left( e^{it/\s_Z} - 1\right) + 
{u^2H''(u)\over 2} \left( e^{it/\s_Z} - 1\right)^2 + R_2(ue^{it/\s_Z},u),
$$
where, according to (2.13), 
$$
| R_2(ue^{it/\s_Z},u)| \le 
h_0(u,r)\cdot  { u^3 | e^{it/\s_Z} - 1|^3 \over r^3}
\cdot {1\over 1-u | e^{it/\s_Z} - 1|/r}.
\eqno (2.14)
$$

It follows from (2.4) and (2.5) that if $k\to\infty$, then  in all of the three 
asymptotic regimes (1.12), (1.14) and (1.22) we have $\s_Z\to \infty$ . Therefore for any given $t$,
$$
| e^{it/\s_Z} - 1|= O(t/\s_Z), \quad (k,x_k)\to\infty.
\eqno (2.15) 
$$
Then we can write that 
$$
H(u e^{it/\s_Z}) - H(u) =
uH'(u) {it \over \s_Z} - 
{t^2\over 2 \s_Z^2} \left( uH'(u) + u^2 H''(u)\right)
+  R_2(ue^{it/\s_Z},u) + S_2(u,t)
\eqno (2.16)
$$
and 
$$
\Phi_Y(t) = \exp\left\{ - {t^2\over 2} + x R_2(ue^{it/\s_Z}, u) + x S_2(u,t)\right\},
\eqno (2.17) 
$$
where 
$$
|xS_2(u,t)|= \left( uH'(u) + u^2H''(u)\right) 
O\left( {t^3\over \s_Z^3}\right) = O(t^3/\s_Z), \quad (k,x_k)\to\infty.
\eqno (2.18)
$$
In (2.17) and (2.18) we have used twice definition of $\s_Z$ (2.5).

Taking into account  elementary upper estimate
$$
h_0(u,r)= \max_{\phi\in [0,2\pi)} | H(u + r e^{i\phi})| 
\le  \max_{\phi\in [0,2\pi)} 
\sum_{k=0}^\infty V_k{ | u+r e^{i\phi}|^k\over k!} \le H(u+r),
$$
and the upper bound (2.13), we can write that  in the case of infinite $u_0=\infty$, $r=1$,
starting form certain $k_0$, 
$$
| xR_2(ue^{it/\s_Z},u)| \le x Q^{(1)}(u) { |t|^3\over \s_Z^3} , \quad 
Q^{(1)}(u) = 2u^3 H(u+1) .
\eqno (2.19)
$$
In the case of finite $u_0$, with the choice of $r= (u_0-u)/2$,
 we can write that 
$$
| xR_2(ue^{it/\s_Z},u)| \le xQ^{(2)}(u)  {|t|^3\over \s_Z^3} ,
\quad Q^{(2)}(u)= 
{16 u^3_0 H((u_0+u)/2))}.
\eqno (2.20)
$$

Taking into account (2.18) together with either (2.19) or (2.20), we  conclude that  two last terms of (2.17)
vanish in the limiting transitions 
$(k,x_k)_\infty \to\infty$ and $(k,x_k)_\chi \to \infty$ and therefore relation (2.9) is true in this two asymptotic regimes.

Regarding the limiting transition $(k,x_k)_0\to\infty$ when $u_k\to \infty$, we 
see that if   (1.6) is true, then 
the right-hand side
of either (2.19) or (2.19) vanishes. This completes the proof of Lemma 2.1. 
Let us note that Lemma 2.1 is also true in the limit $(k,x_k)_0\to\infty$,  where the sequence $x_k\to 0$  in a way
 that (1.11) is verified. 
 \hfill  $\Box$

 \vs
 \subsection{Local Limit Theorem}
 
 Let us show that the random variables $Z^{(x,u)}$ verify the Local Limit Theorem. 
 \vs 
 {\bf Lemma 2.2.}
 {\it  
Given a sequence of $(x_k,\chi_k)_{k\in \bN}$, we consider $u_k$ such that (1.9) is verified for all $k\in \bN$. Then relation 
$$
P(Z^{(x_k,u_k)}= k) = {1\over \sqrt{2\pi } \s_Z}\big(1+o(1)\big), \quad  (k,x_k)\to\infty
\eqno(2.21)
$$
holds in the  asymptotic regimes
$(k,x_k)_\infty\to\infty$ (1.12)  and $(k,x_k)_\chi\to\infty$ (1.14).
If $H(u)$ is such that conditions (1.6), (1.7) and (1.8) are verified, then (2.21) is also true in the limiting transition $(k,x_k)_0\to\infty$ (1.22).
}

\vs 
{\it Proof of Lemma 2.2.}
We combine  relation  (2.18) with 
  arguments developed by   T. Tao in \cite{Tao} for  the proof of the Local Limit Theorem
for sums of independent random variables.
Taking  mathematical expectation of both parts of equality
$$
{\bf I}_{Z=k}(\omega) = {1\over 2\pi} 
\int_{-\pi}^\pi e^{iyZ} e^{-iyk} dy, 
$$ 
we get by the Fubini's theorem that 
$$
P(Z=k) = {1\over 2\pi}\int_{-\pi}^\pi \E \left( e^{iy(Z-\E Z)}\right) \ e^{-iy(k-\E Z)} \, dy
$$
$$
= {1\over 2\pi \s_Z} 
\int_{- \pi \s_Z}^{\pi \s_Z} \Phi_Y(t) \, 
e^{it(k-\E Z)/\s_Z} \, dt,
\eqno (2.22)
$$
where $Y= Y^{(x,u)}$
is determined by (2.8). 
The key point of the proof is to show that 
the following difference vanishes,
$$
D_k= \int_{- \pi \s_Z}^{\pi \s_Z} \Phi_Y(t) \, 
e^{it(k-\E Z)/\s_Z} \, dt
-
\int_{-\pi \s_Z}^{\pi \s_Z} e^{it(k-\E Z)/\s_Z - t^2/2}\, dt = o(1)
\eqno (2.23)
$$
as $k$ tends to infinity. 
Then (2.21) will follow from this convergence, classical identity
$$
{1\over 2\pi} \int_{-\infty}^{\infty} e^{it(k-\E Z)/\s_Z - t^2/2}\, dt = {1\over \sqrt{2\pi}}
e^{-(k - \E Z)^2/(2\s_Z^2)},
$$
and obvious estimate
$$
\left| \int_{| t|> \pi \s_Z} e^{-t^2/2 + i \a t} dt  \right| = o(1), \quad \s_Z\to\infty.
$$

We split $D_k$ into two parts
and consider first the difference
$$
D_k^{(1)}= \int_{\vert t \vert \le y_k} \left( \Phi_Y(t) - 
e^{-t^2/2}\right) e^{-it(k-\E Z)/\s_Z}\, dt.
$$
It follows from (2.17), (2.19) and (2.20)  that  
$$
\vert D_k^{(1)}\vert  \le 
\int_{\vert t \vert \le y_k} e^{-t^2/2} \vert \exp
\left\{ x_kR_2 + x_kS_2 \right\}- 1\vert  
\, dt 
$$
$$
=O\left( x_k Q^{(l)}(u_k) {y_k^k\over \s_Z^3}\right)+  O\left({x_k y_k^4\over \s_Z}\right) = o(1), \quad (k,x_k)\to\infty.
$$
It is clear that with the choice of  $y_k= \s_Z^{1/8}$  we have 
$$
|D^{(1)}_k| = o (1)
\eqno (2.24) 
$$
in the limiting transitions $(k,x_k)_\infty\to \infty$ and $(k,x_k)_\chi\to\infty$. 

Regarding
the second part of the integral (2.23), we can write that
$$
\vert D_k^{(2)}\vert \le  
\int_{y_k< \vert t\vert  \le \pi \s_Z}
e^{-t^2/2}\,  dt + \int_{y_k< \vert t\vert  \le \pi \s_Z}
\vert \Phi_Y(t)\vert \, dt. 
\eqno (2.25)
$$

\vs
Let us consider the last integral  in  the limiting transitions $(k,x_k)_\infty\to\infty$ (1.12) 
and $(k,x_k)_\chi\to\infty$ (1.14) when the parameter $u_k$ remains finite as $k\to\infty$. 
Using  definitions (1.5) and (2.1) together with  relations (2.2) and (2.11), we can  write that 
$$
\vert \Phi_Y(t) \vert = \left\vert \exp
\left\{ x_k \sum
{V_j } u_k^j 
\big(  \cos(jt/\s_Z) -1)
+ i \sin(jt/\s_Z)\big)/j! \right\}\right\vert 
$$
$$
= \exp\left\{ x_k \sum_{j=1}^2 u_k^j {V_j} 
\left( \cos\left({jt/ \s_Z}\right) -1 \right) /j! 
+ x_k T_3(x_k,u_k,t)\right\},
\eqno (2.26)
$$
where the remaining term
$$
T_3(x_k,u_k,t)= \sum_{j=3}^\infty V_j u_k^j 
\left(  \cos\left(jt/ \s_Z\right) -1\right)/j!
$$
is negative,  $T_3(x_k,u_k,t)\le 0$ for all $t\in [-\pi \s_Z, \pi\s_Z]$.

Taking into account elementary inequalities
$$
\cos\left( {j t\over \s_Z}\right) \le 1 - {j^2 t^2\over 24 \s_Z^2} \le 
1 - {j ^2y_k^2\over 24 \s_Z^2}, \quad \forall t: 
\ y_k\le  \vert t\vert \le \pi \s_Z, \quad j = 1, 2, 
\eqno (2.27)
$$
we deduce from (2.26) that
$$
\Delta_k =  \int_{y_k< \vert t\vert  \le \pi \s_Z}
\vert \Phi_Y(t)\vert \, dt \le 2\pi \s_Z
\exp\left\{ - { x_k u_k  (V_1+ 4u_kV_2) y_k^2\over 24 
\s_Z^2}
\right\}.
$$
Remembering definition (2.5) of $\s_Z^2$ and the value $y_k = \s_Z^{1/8}$, we conclude that 
$$
\Delta_k  \le 2\pi \s_Z
\exp\left\{ - {  V_1+ 4u_kV_2\over 24 (H'(u_k) + u_k H''(u_k))} 
\s_Z^{1/4}
\right\} = 2\pi \s_Z \exp\{ - A_k \s_Z^{1/4}\}.
$$
It is clear that in the limiting transition $(k,x_k)_\chi\to \infty$ when  $u_k$ converges as $k\to\infty$ to a finite solution  $u>0$ of equation (1.17),
 we get $\Delta_n = o(1)$ .  
  
Let us consider the limiting transition  $(k,x_k)_\infty \to \infty$ (1.12). In this case $\chi_k\to \infty$ and therefore  $u_k\to 0$. If 
$V_1 \neq 0$, then $\lim_{n\to\infty} A_k = 1/24$ and if $V_1 = 0$, then
$\lim_{k\to\infty } A_k = 1/6$. Then 
$\Delta_k = o(1), \ (k,x_k)_\infty \to \infty.
$
Returning to (2.25), we conclude that in these two limiting transitions,
$$
\vert D_k^{(2)}\vert = o(1), \quad (k,x_k) \to\infty.
$$
Combining this relation with (2.24), we get (2.23).  Then (2.22) implies (2.21).  The first part of Lemma 2.2 is proved.

Now let us consider the  limiting transition $(k,x_k)_\infty\to\infty $ (1.22)  when $u_k\to\infty$.
It follows from (2.17), (2.19) and (2.20) that 
$$
|D_N^{(1)}| \le \int_{|t|\le y_N} | e^{xR_2 + x S_2} - 1| dt
$$
$$
= \int_{|t|\le y_N} | \big(xR_2 + x S_2\big)(1+o(1))| dt
\le y_N^4 {\CT}_N(x,u),
$$
where 
$$
Q_H (x,u) =  \begin{cases}
{xu^3 H(u+1)/\s_Z^3 + 1/\s_Z}
   , & \text{if $u_0=\infty$ } , \\
xu^3 H((u_0+u)/2)/((u_0-u)^3\s_Z^3 + 1/\s_Z, & \text {
if  $u_0<\infty $.} 
\end{cases}
\eqno (2.28) 
$$

To estimate $D_N^{(2)}$, we start with  the first 
integral of the right-hand side of (2.25). 
Remembering that 
$$
H(ue^{i\a}) = \bE \exp\{u(\cos \a + i \sin \a) W\}
=
\bE e^{uW \cos \a} \big( \cos (uW \sin \a) + i \sin(uW \sin \a)\big),
$$
we write that
$$
| \exp\left\{x H(u e^{i\a}) - x H(u) \right\}|\le 
\exp\left\{ x\bE \left( e^{uW\cos \a}- e^{uW}\right)\right\},
$$
where $\a= t/\s_Z$. 

Using the  upper bound  $\cos s\le 1 - s^2/12$, $|s|\le \pi$, we can write inequality
$$
|D_N^{(2)}| \le 2\pi \s_Z \exp\left\{ x \bE e^{uW} 
\left(e^{-uW y^2/(12\s_Z^2)} - 1 \right) \right\}.
$$
Then
$$
| D_N^{(2)}|  \le 2\pi \s_Z \exp\left\{ x \bE e^{uW} 
\left(e^{-uW y^2/(12\s_Z^2)} - 1 \right){\hbox{I}}_{\{W\ge w_0\}} \right\}
$$
$$
\le 2\pi \s_Z \exp\left\{ x\left( e^{-w_0 uy^2/(12\s_Z^2)} - 1\right) \beta H(u)\right\}
\le 2\pi \exp\left\{ - {x u w_0 \beta H(u)\over 24  \s_Z^2} y^2 + \ln \s_Z\right\}.
\eqno (2.29)
$$
In the last estimate, we have used elementary inequality
$$
e^{-s}-1\le -{s\over 2}
$$
that is true for sufficiently small $s>0$.

Gathering (2.28) and (2.29), we see that if 
the value of $y_k$ infinitely increases as $k\to\infty$ and and verifies  the following chain of asymptotic estimates,
$$
{ \s_Z^4\over (x_k u_k H(u_k))^2} \big( \ln \s_Z\big)^2 \ll 
y_k^4 \ll {1\over Q_H(x_k,u_k)}\ , \quad (k,x_k)_0\to\infty,
\eqno (2.30)
$$ 
then (2.23) holds. 
It is not hard to see that conditions (1.7), (1.8) are sufficient to have 
the left-hand side of (2.30) much less than the right-hand side of (2.30). 
Lemma 2.2 is proved. $\Box$

\subsection{Proof of Theorem 1.1 and Corollaries A and B} 

Taking into account relation   (2.21) of Lemma 2.2 
and using definitions (2.1), (2.2) and (2.7),
we obtain the following   asymptotic relation,
$$
M_k(x)=
 {k!\over u^k}\cdot { \exp\left\{ {x_k(H(u_k)-1)} \right\}
 \,  \over \sqrt{2\pi } \s_{Z}^{(k)}} \big(1+o(1)\big)
 \,  , \quad (k,x_k)\to\infty,
\eqno (2.31)
$$
where where  
$\s^{(k)}_Z = \sqrt{ x_k u_k \left( H'(u_k) + u_k H''(u_k)\right)}$ and  
$u_k$ is determined by (1.9). 
Rewriting  equality (1.9) in the form
$$
x_k= {k\over u_kH'(u_k)},
$$
we deduce from (2.31) that 
$$
M_k(x_k) =  
 {k!\over  \sqrt{2\pi }\s_Z^{(k)} u_k^k}\, \exp\left\{ k{ 
 H(u_k)-1\over u_kH'(u_k)}\right\} (1+o(1))
, \quad (k,x_k)\to \infty. 
\eqno (2.32)
$$
Using the Stirling formula
$$
k! = \sqrt {2\pi k} \left( {k\over e} \right)^k (1+o(1)),\quad k\to\infty
\eqno (2.33)
$$
we  transform relation (2.32) into the following asymptotic equality
$$
M_k(x_k)= \left( { k\over  u_k }\right)^k \exp
\left\{ k\left({H(u_k)-1\over u_k H'(u_k)}-1\right)\right\}   \sqrt{ {k\over x_k u_k(H'(u_k) + u_kH''(u_k))}}(1+o(1)), 
\eqno (2.34)
$$
that is valid in all fo the three asymptotic regimes (1.12), (1.14) and (1.22). Then (1.10) obviously follows from (2.34). {Theorem 1.1} is proved. 
\hfill $\Box$

\vs  
Let us prove two corollaries of Theorem 1.1 given by (1.13) and (1.15) in two asymptotic regimes (1.12) and (1.14), respectively.

A. Consider the case when 
$1\ll k\ll x_k$.  Then  $\chi_k\to \infty $ as $k\to\infty$ and it follows from (1.9) that 
$u_k $ converges to zero as $k\to\infty$.  Then we can write that 
$$
H(u_k)-1 - u_kH'(u_k) = -{ u^2_k V_2\over 2} + o(u^2_k), \quad (k,x_k)_\infty\to\infty.
$$
If $V_1\neq 0$, then $H'(u_k)= V_1(1+o(1))$ and 
$$
\exp\left({H(u_k)-1\over u_k H'(u_k)}-1\right) =1+ o(1), 
\quad (k,x_k)_\infty\to\infty.
$$
Remembering (1.9), we see that $u_k= k(x_kV_1)^{-1} (1+o(1))$. Then we deduce from (2.34) relation
$$
M_k(x_k) = \big(x_k V_1 (1+o(1))\big)^k, 
\quad (k,x_k)_\infty\to\infty.
\eqno (2.35) 
$$
If $V_1=0$, then $H'(u_k) = V_2 u_k(1+o(1))$
and relation (1.9) implies equality
$$
u_k= \sqrt{k\over x_k V_2} (1+o(1)), \quad (k,x_k)_\infty\to\infty.
$$
Taking into account that
$$
\exp\left({H(u_k)-1\over u_k H'(u_k)}-1\right) = e^{-1/2}(1+o(1)), 
\quad (k,x_k)_\infty\to\infty,
$$
we deduce from (2.34) asymptotic equality
$$
M_k(x)= \left({ \sqrt{{ x k V_2\over  e}} }\big(1+o(1)\big)\right)^{k}, \quad 
(k,x_k)_\infty\to\infty.
\eqno (2.36) 
$$
Relations (2.35) and (2.36) give (1.13).

\vs

B. If $x_k/k\to  \chi>0$, then $u_k\to u$
with $u$ determined by the Lambert-type equation (1.9). 
It follows from (2.34)  that in this case 
$$
M_k(x_k) =   \left( x_kH'(u)  \exp\left\{ {H(u)-1\over u H'(u)} 
-1 \right\}\big(1+o(1)\big) \right)^k , \ \ (k,x_k)_\chi\to \infty. 
\eqno (2.37)
$$
Then (1.15) and (1.16) follow from (2.37).

Let us note that if we determine $x_k$ by relation $x_k= k\chi$, 
then (2.34) implies the following asymptotic equality
$$
M_k(x_k) =  
{1\over \sqrt{1 + \chi u^2 H''(u)}} \left( {k\over u } \exp\left\{ {H(u)-1\over u H'(u)} 
-1 \right\} \right)^k(1+o(1)) , \ \ (x,k)_\chi\to \infty 
\eqno (2.38)
$$
that is 
 more informative  than (2.37).


\section{The third asymptotic regime} 

In this section we consider  asymptotic behavior of the moments 
$M_k(x_k)$ in  the limit \mbox{$(k,x_k)_0\to \infty$} (1.22) when
$x_k$ is much smaller that $k$. 
We concentrate ourself on several important particular cases of the probability
distribution of the weights $X_j$.

\subsection{Gamma distribution}

Assuming that $X_j$ follows the Gamma distribution with density
 $$
 f^{(m,\theta)} (x) = {x^{m-1}\,  e^{-x/\theta} \over \th^m \Gamma(m)},
 \quad x>0, \ m>0, \ \th>0,
 $$
 where 
 $\Gamma(m) = \int_0^\infty x^{m-1} e^{-x} dx$ \cite{KPST}, 
 it is not hard to see  that (1.7) holds and that 
 $$
 H(u)= {1\over (1 - \th u)^{m}}, \quad 0\le u < {1/ \th}.
 \eqno (3.9)
 $$
 We denote by $M_{k}^{(\Gamma)}(x)$ the moments (1.5) 
 with the weights $V_j= \int_0^\infty x^j f^{(m, \theta)}(x) dx$.  
 Let us consider random variables $Z^{(x,u)}$ such that   
 $$
 P(Z^{(x,u)} =k) = M_k^{(\Gamma)}(x) {u^k\over k! \, G(x,u)} =  S_k(x) {u^k\over G(x,u)}
 \eqno(3.10)
 $$
with 
$$
G(x,u) = \exp\left\{ x(H(u)-1)\right\} = 
\exp\left\{ x\left( (1-\theta u)^{-m}-1\right)\right\}
$$   
and 
$$
\E Z = {m\theta xu\over (1- \theta u)^{m+1}}
\quad \hbox{and} \quad 
\s_Z^2={m\theta xu\over (1- \theta u)^{m+1}}
+ 
xu^2 {m(m+1) \theta^2\over (1- \theta u)^{m+2}}.  
$$
 
 The Lambert-type equation (1.9) takes the form
 $$
 {m\th u\over (1-\th u)^{m+1}} = {k\over x}\, .
 \eqno (3.3)
 $$

Assume that conditions (1.6), (1.7) and (1.8) are verified. 
 Then Lemma 2.2 is true in the asymptotic regime $(k,x_k)_0\to \infty$ (1.22).
 In this case the solution of (3.3) admits the following  expansion,
  $$
u_k = {1\over \th } + {1\over \th} \left({ m x_k\over  k}\right) ^{1/(m+1)}\big(1+o(1)\big) ,
$$

Denoting  by $M_{k}^{(\Gamma)}(x)$ the moments (1.2) 
 with the weights $V_j= \int_0^\infty x^j f^{(m,\th)}(x) dx$, we can write that 
 $$
 M_k^{(\Gamma)}(x_k)=\left({ mx_k\over k}\right) ^{1/(m+1)}
 \left( { k\th \over e 
 \left(
 1 +  \left({ mx_k/ k}\right) ^{1/(m+1)}
 (1+o(1))\right)}
 \right)^k \ (1+o(1))
 $$
 $$
 =\left( {k\th \over e }\big(1+o(1)\big)\right)^k, \quad (k,x_k)_0\to\infty.
\eqno(3.4)
$$ 
 This proves proposition (1.23).

   Let us show that  conditions (1.6), (1.7) and (1.8) are verified when  the Gamme distribution is such that $m>8$.
   We start with (1.8). We can write  that
$$
{1\over \th^m \G(m)}\int_{x\ge  w_0} e^{ux} x^{m-1} e^{-x/\th } dx 
$$
$$
= 
{1\over (1/\th-u)^m}\left( 1 - 
{1\over \th^m  \G(m)} \int_{s< w_0(1/\th -u)} s^{m-1}e^{-s} ds\right)
$$
$$
 \ge {1\over (1/\th -u)^m}\left( 1 -  
{1\over \th^m\G(m)}\int_{s< w_0} s^{m-1}e^{-s} ds\right)
$$
$$
\ge 
{1\over (1/\th -u)^m}\left( 1 -  
{1\over \th^m\G(m)}\int_{s< w_0} s^{m-1} ds\right).
$$
We conclude that if $w_0$ is such that
$$
{w_0^m\over m \th ^m  \G(m)}<{1\over 2},
$$
then  relation (1.8) holds with $\beta = 1/2$.

 Instead of (1.6), let us check whether relations  (2.30) are satisfied. 
This happens  when
$$
\s_Z \left( \ln \s_Z\right)^2 \left( 16x H\left({u+u_0\over 2}\right) + (u_0-u)^3 \s_Z^2\right)\ll (u_0-u)^3 \big(x uH(u)\big)^2
\eqno (3.5)
$$
holds in the limit $u\to u_0$. 
Using equality $H((u+u_0)/2))= 2^m(1-\th u)^{-m}$ and 
observing that 
$$
\s_z^2= = {x m (m+1) \th^2\over (1-\th u)^{m+2}}\big(1+ o(1)\big), \quad u\to 1/\th,
$$
we conclude that (3.5) is verified when the following relation holds,
$$
{2^{m+4} x^{3/2} \sqrt{m(m+1)} (m+2)^2 \th \over (1-\th u)^{3m/2 + 1}}
\left(\ln(1-\th u)\right)^2\ll{x^2\over (1-\th u)^{2m-3}}, \quad u\to 1/\th .
$$ 
This is true 
under condition 
$m>8$. 
It is interesting to note that condition (2.30) imposes restriction $m>8$ on $m$ only while 
the leading term of the asymptotic expression (3.4) depends on $\th$ only. 

\subsection{Normal (Gaussian) distribution}
Let us consider first the centered random variables. 
If  $W_j\sim {\CN}(0,V_2)$, 
then 
 $$
H(u)= e^{ {V_2 u^2/ 2}}, \quad u\in \bR.
 \eqno (3.6)
$$
In this case equation (1.9) takes the form of Lambert equation (1.19)
$$
se^{s} = t={1\over \chi},
\eqno (3.7)
$$
where we denoted $s= u^2 V_2 /2$. 
In the third asymptotic regime (1.22) when $x_k= o(k)$, we have $\chi= \chi_k\to 0 $, $k\to\infty$.
It is known that the asymptotic expansion of $s=s(t)$  is as follows \cite{dB},
$$
s(t)= \ln(t) - \ln \ln(t)(1+o(1)), \quad t\to\infty.
\eqno (3.8) 
$$ 

 It is not hard to show that conditions (1.6), (1.7) and (1.8) are verified
 by $H(u)$ (3.6). 
Then 
$$
u_k= \sqrt{{2(\ln (k/2) - \ln  x_k)\over V_2}}\big(1+o(1)\big).
\eqno (3.9)
$$
Taking into account that 
$$
{H(u_k)-1\over u_kH'(u_k)}\to 0, \quad u_k\to\infty,
$$
we deduce from (1.10) with the help of (3.9) that 
 $$
 M_{2k}^{(\CN(0,V_2))} (x_k)= \left( { 2k\over e u}(1+o(1))\right)^{2k}= 
 \left( {\sqrt { 2V_2} k\over e \sqrt{\ln (2k) - \ln x_k}} (1+o(1))\right)^{2k}, \quad (k,x_k)_0\to\infty.
 \eqno (3.10)
 $$  Relation (1.25) is proved. 
 
 \vs
  
 Now we consider the general case $X_j \sim \CN(a,\sigma^2)$. It is easy to see that conditions (.16), (1.7) and (1.8) are verified. Then  
 $$
 H(u) = e^{ua + u^2 \s^2/2}, \quad u\in\bR
 $$ and 
 equation (1.9)
 takes the following form slightly different from (3.7), 
 $$
 (au+ \s^2 u^2) e^{au + u^2\s^2/2} =  t, \quad t={1\over \chi}.
 \eqno (3.11)
 $$
 To study asymptotic expansion of $u$, we denote the left-hand side of (3.11) by $\tilde F(u)$ and use the denotation 
 $\tilde u= \tilde u(t)$ for the solution of the corresponding equation $\tilde F(u)= t$. 
 
 We introduce two auxiliary functions, $\hat F(u)$ and $\check F(u)$ by the following formulas, 
 $$
\hat F(u) =  (au + \s^2 u^2/2) e^{au + u^2\s^2/2}< \tilde F(u) < \check  F (u) = (2au + \s^2 u^2) e^{au + u^2\s^2/2} .
 $$
 Then clearly 
 $$
 \check u(t) \le \tilde u(t) \le \hat u(t),
 $$
 where $\check u$ and $\hat u$ denote the solutions of equations $\check F(u)=t$ and $\hat F(u)=t$, respectively. 
 
 Remembering asymptotic  expansion (3.8),  
 we conclude that $\check u(t)$ and $\hat u(t)$ are such that
 $$
 a\hat u + \s^2 \hat u^2/2 = \ln t (1+o(1)) \quad \hbox{and} \quad a\check u +\s^2 \check u^2/2 = \ln(t/2)(1+o(1)).
 $$
Then
 $$
 {\sqrt {2\ln(t/2)}\over \s} (1+o(1)) \le \tilde u(t) \le {\sqrt {2  \ln t}\over \s} (1+(1))
 $$
and we can write that
 $$
 \tilde u (t) = {\sqrt {2 \ln t}\over \s}\big( 1+o(1)\big)= {\sqrt {2  \ln(k/x)}\over \s} \big(1+o(1)\big).
 $$
 Substituting this expression into (1.10) (see also (2.5)),
 we get relations
 $$
 M_k^{(\CN(a,\s^2))}(x_k) = \left({ k\s\over e \sqrt {2 (\ln k - \ln x_k)}} \big(1+o(1)\big)\right)^k, \quad (k,x_k)_0\to\infty.
 \eqno (3.12)
 $$ 
 Relation (1.24) is proved. 
 Let us  note that $a$ does not alter the leading term of $M_k(x_k)$. 
 Moreover, we get the same asymptotic behavior as in the case of centered Gaussian distribution (3.10).

\subsection{Centered Bernoulli  and triangular distributions}

Let us consider the case when independent random variables $X_j$ (1.1) follow 
the centered Bernoulli probability distribution,
$$
X_j=
 \begin{cases}
  1 , & \text{with probability $1/2$} , \\
-1, & \text {with probability $1/2\, $} 
\end{cases}, \quad j\in \bN. 
\eqno (3.13)
$$
In this case, the moments (1.2) take the form 
$$
 M^{(\CB_0)}_{2k} (x) = 
 \sum_{(l_2, l_4 , \dots, l_{2k})^*}^{2k} 
 x^{l_2 + l_4 + \dots + l_{2k}}   \hat B_{2k}(l_2, l_4, \dots, l_{2k}),
 $$
 where 
 $$
 \hat B_{2k} (l_2, l_4 , \dots, l_{2k})= { (2k)!\over (2!)^{l_2}l_2!\,  
 (4!)^{l_4} l_4! \cdots ((2k)!)^{l_{2k}} l_{2k}!}
 $$
 and the sum runs over $l_{2i}\ge 0$ such that $2l_2 + 4l_4 + \dots + 2kl_{2k}= 2k$. 
 Then 
 $$
 M_{2k}^{(\CB_0)}(1) = \hat B^{(2k)} = 
 \sum_{(l_2, l_4 , \dots, l_{2k})^*}^{2k} 
 \hat B_{2k}(l_2, l_4, \dots, l_{2k}),
 \eqno (3.14)
 $$
 represents the number of all possible partitions of a set of $2k$
 elements into subsets  of even cardinality \cite{OEIS}. One can refer to 
 $\hat B^{(2k)}$ as to the  Bell numbers of even partitions.

Definition  (3.13) implies that  
 $$
 H(u) = \bE e^{uX}= \ch(u)
 $$ 
and the Lambert-type equation (1.9) takes the form 
$$
u \sh(u) = t, \quad t=  {2\over \chi}.
\eqno (3.15) 
$$
Taking into account elementary bounds 
$$
{e^u\over 4} \le \sh (u) \le {e^u\over 2}, \quad u>1,
$$
it is easy to show with the help of the arguments of the previous subsection that (3.15) implies 
asymptotic equality
$$
u(t) = \ln t (1+o(1)), \quad t\to\infty
$$
It not hard to see that conditions (1.6), (1.7) and (1.8) are verified in the case of centered Bernoulli distribution (3.13). Then it  follows from  Theorem 1.1 that 
 $$
  M_{2k}^{(\CB_0)}(x_k) = \left( { 2k\over e(\ln(2k) - \ln x_k)}(1+o(1))\right)^{2k}, \quad (k,x_k)_0\to\infty.
  \eqno (3.16)
  $$
  Relation (3.16) implies (1.26).
Relation (1.26), when taken with $x=1$, coincides with the first terms of the asymptotic expression 
 for Bell numbers $B^{(2k)}$ (1.27) in the limit $k\to\infty$,  known since 50-s \cite{MoWy}. 
 
\vs
Let us consider the case when $X_j$ have triangular distribution $X_j\sim T(\nu)$ with the density of the form
$$
f^{(\nu)}(x) = {1\over \nu} 
 \begin{cases}
  1+x/\nu , & \text{if  $x \in (-\nu,0]$} , \\
1-x/\nu, & \text {if  $x \in (0,\nu] $}\\
0, & \text{otherwise} 
\end{cases}.
$$
It is easy to see that
$$
V_{2i}^{(\nu)}= \int_{-\nu}^\nu x^{2i} f^{(\nu)}(x) dx  = 
{\nu^{2i}\over (i+1)(2i+1)} = \nu^{2i} V_{2i}^{(1)}.
$$
It follows from definition (1.2) that 
$$
M_{2k}^{(T(\nu))}(x)= \nu^{2k} M_{2k}^{(T(1))}.
$$

Regarding the case of $\nu=1$, we can write that  
$$
H(u)= {2\over u^2}\big( \ch(u) - 1\big)
$$
and equation (1.9)  takes the form
$$
\tilde F(u) = t, \quad t= {2k\over x}= {2\over \chi},
\eqno (3.17)
$$
where we denoted 
$$
\tilde F(u) = {2\over u} \sh(u)  - {4\over u^2} \big( \ch(u) - 1\big).
\eqno (3.18)
$$

\vs
In Section 4, we show that solution $\tilde u$ of (3.17) admits the following asymptotic expansion in the limit $\chi\to 0$, 
$$
\tilde u = \ln(k/x)(1+o(1)).
\eqno (3.19) 
$$
Then 
$$
M_{2k}^{(T(1))} (x_k)= \left( {2k\over e(\ln(2 k) - \ln x_k)}\big( 1+ o(1)\big) \right)^{2k}, \quad (k,x_k)_0\to\infty.
$$
that coincides with (1.26) and (3.16). This coincidence is similar to that observed in the normal distribution, the right-hand side of (3.10)
is the same of (3.12) with $k$ replace by $2k$. 

Finally, we get for the the moments of compound Poisson distribution with 
$X_j\sim T(\nu)$ the following asymptotic expression,
$$
M_{2k}^{(T(\nu ))} (x_k)= \left( {2k\nu \over e(\ln(2 k) - \ln x_k)}\big( 1+ o(1)\big) \right)^{2k}, \quad (k,x_k)_0\to\infty.
\eqno (3.20)
$$
 
\vs



\section{Auxiliary statements and supplementary results}

In this section we collect auxiliary facts used in the proof and discuss  a number of additional 
results related with Theorem 1.1.

\subsection{Taylor expansion}

In this subsection we prove relations (2.12) and (2.13).
It follows from (1.5) that the series 
$$
H(z) = \sum_{k=0}^\infty V_k {z^k\over k!} 
$$ 
converges for any $z$ from the open ball $B(0, u_0)$ and by definition is analytic in $B(0, u_0)$. It is (infinitely) holomorphic
and the Cauchy integral formula is true,
$$
H^{(k)}(u) = {k!\over 2\pi i}  
\int_C {H(s)\over (s-u) ^{k+1}}ds , \quad s= u + r e^{i\phi}, 
\quad r\le {u_0 - u\over 2}. 
\eqno (4.1)
$$ 
Therefore we can write  that  for any $z$ and $u$ from $ B(0, u_0)$, 
$$
H(z) = \sum_{j=0}^\infty H^{(j)} { (z-u)^j\over j!} = 
\sum_{j=0}^2 H^{(j)} { (z-u)^j\over j!} + R_2(z,u),
$$
where 
$$
R_2 (z,u) = \sum_{j=3}^\infty H^{(j)} { (z-u)^j\over j!}.
$$
It follows from (4.1) that 
$$
\vert H^{(j)} (u) \vert \le {j! h_0(u,r)\over r^k}, 
\quad h_0(u,r) = \max_{s: \, \vert s-u \vert = r} \vert H(s)\vert.
$$
Finally, we obtain the following estimate,  
$$
\vert R_2(z,u)\vert \le \sum_{j=3}^\infty { h_0(u,r) \vert (z-u)\vert ^j
\over r^j} = h_0(u,r) { \vert z-u\vert ^3\over r^3 \left( 1 - 
\vert z-u\vert /r\right)}.
$$
This proves relations (2.12) and (2.13). 

\subsection{Lambert-type equation for triangular distribution}

It follows from the definition (3.18) that
$$
\tilde F(u) = \sum_{l=1}^\infty {4l  \over (2l+2)!}\, u^{2l}.
$$
Then we can write inequality 
$$
\tilde F(u) \ge 4\sum_{l= 1}^\infty {{u^{2l}\over (2l+2)!}} = 
{4\over u^2}\big(\ch(u) - 1\big) - 2.
\eqno (4.2)
$$
Combining (3.18) with (4.2), we conclude that  
$$
{4\over u^2}\big(\ch(u) - 1\big) \le {\sh(u)\over u} + 1
$$
and 
that
$$
\tilde F(u) \ge {\sh(u)\over u} - 1= \hat F(u) .
$$
Denoting by $\hat u = \hat u (t)$  solution of equation 
$$
\hat F(u)= t,
\eqno (4.3)
$$
we conclude that $\tilde u(t) \le \hat  u(t)$, where $\tilde u(t)$ is the solution of equation (3.17).

Taking into account that for  $u\ge \ln 2$
we have
$
\sh(u) \ge e^u/4$,
we can write that 
$\hat u(t) \le \bar u(t)$, where $\bar u(t)$ is a solution of equation 
$$
{e^u\over u}= 4 t+4.
\eqno (4.4)
$$
It follows from (4.4) that 
$$
u= \ln (4t+4)  + \ln u \le \ln (4t+4) + u
/2.
$$
The last inequality is true for sufficiently large $t$ because $u \ge \ln (4t+4)$. 
Then 
$
 u\le 2 \ln (4t+4) 
 $ and 
$$ 
\ln u \le \ln \ln (4t+4) + \ln 2.
 $$
 Then we get the upper bound 
 $$
 \tilde u(t)\le \bar u(t) \le \ln (4t+4) + \ln \ln (4t+4) + \ln 2.
 \eqno (4.5) 
 $$

Let us find the lower bound for $\tilde u(t)$. 
Remembering definition (3.18), we can write that 
$$
\tilde F(t) \le \check F(t)= {2\over u} \sh(u)
$$
and therefore $\tilde u(t) \ge \check u(t)$, where $\check u(t)$ is a solution of equation
$$
\check F(u) = t.
$$

Taking into account inequality $2\sh(u) \le e^u$,
we can write that $\check u(t) \ge \dot u(t)$, where 
$\dot u(t)$ is a solution of equation 
$$
{e^u\over u} = t.
$$
It is clear that $\dot u = \ln t + \ln{\dot u}$ and therefore $\dot u(t) \ge \ln t$.
Then 
$
\tilde u(t) \ge \ln t.
$
This relation together with (4.5) shows that 
$$
\tilde u(t) = \ln (t) \big(1+o(1)\big), \quad t\to \infty.
\eqno (4.6) 
$$
This relation implies (3.19).

\subsection{Exponential distributions and combinatorial polynomials}
If $X_j$ follow the exponential distribution, $X_j\sim \CE(1)$, then $V_k = k!$, $k\in \bN$
and 
$$
 H(u)= \sum_{k=0}^\infty V_k {u^k \over k!} = {1\over 1-u}, \quad u\in [0, 1).
 \eqno (4.7)
  $$
It follows from (1.2) that 
$$
M_k(x) = k!\,  S_k(x),
\eqno (4.8) 
$$
where 
$$
S_k(x)= \sum_{(l_1, \dots, l_k)^*}^{  k} {x^{l_1+\dots + l_k}\over l_1! \dots l_k!}.
\eqno (4.9)
$$

We deduce  from Theorem 1.1 that if $k\to\infty$ and $x= \chi k$, then 
$$
 S_k(x) = {1\over \sqrt{2\pi k}} \sqrt{ {1-u\over 1+u}}\left( {e\over ue^u}
 \right)^k \big( 1+ o(1)\big), \quad k\to\infty,
 \eqno (4.10)
 $$
where $u$ is a solution of equation 
$$
{u\over (1-u)^2} = {1\over \chi},
$$
and therefore
$$
 u = { 2 + \chi - \sqrt { \chi (4+\chi)} \over 2}.
 $$ 

If $k= o(x)$ and $\chi\to\infty$, then 
 $u={k/x} + 3/2 + o(1)$ and   
 $$
 S_k(x) = {1\over \sqrt{2\pi k}} 
 \left( {xe\over k} \big( 1+ o(1)\big) \right)^{k}, \quad 1\ll k \ll x.
 \eqno (4.11)
 $$
 
 \vs 
Using combinatorial    identity (\cite{Rio}, p.183),
$$
 \sum_{(l_1, \dots, l_p) = p} {p!\over l_1! \dots l_p!} =
  {k-1 \choose p-1},
  \eqno (4.12)
 $$
it is easy  to get the following expression for the polynomials $S_k(x)$,
 $$
 S_k(x) =  
\sum_{p=1}^k {x^p\over p!} {k-1 \choose p-1}.
\eqno (4.13)
$$ 
Thus, relations (4.10) and (4.12) determine
 asymptotic behavior of combinatorial polynomials (4.13).

 Let us consider
  $M_k(x)$  (1.2) with  
 $
 V_k = (k-1)!$
 for all $ k\in \bN. 
 $
 We can write that in this case 
 $
 M_k(x) = k!\,  T_k(x)$,
 where 
 $$ 
  T_k(x) =   \sum_{(l_1, l_2, \dots, l_k)^*}^k\,  
 \prod_{i=1}^k {x^{l_i} \over i^{l_i} l_i!}.
$$
We have 
 $$
 H(u)= \sum_{k=0}^\infty V_k {u^k\over k!} = 1 - \ln (1-u).
 \eqno (4.14)
 $$
Introducing auxiliary random variables $Z^{(x,u)}$ such that
$
P(Z^{(x,u)} = k) = T_k(x) {u^k/ G(x,u)}
$
with 
$$
 G(x,u) = \exp\left\{ x(H(u)-1)\right\}  = {1\over (1-u)^x}, \ u\in [0,1),
 $$
 and 
 $$
 \E Z ^{(x,u)} = 
 {xu\over 1-u} \quad \hbox{and} \quad \s_Z^2 = {xu\over 1-u} + {xu^2\over (1-u)^2},
 $$
 we deduce from  Theorem 1.1 
 that if $x_k = \chi k$, then 
 $$
 T_k(x_k) = {1\over \sqrt{2\pi} \s_Z} {G(x_k,u)\over u^k} (1+o(1)), \quad (k,x_k)_\chi\to\infty,
 $$
 where $u$ is determined by (1.9) with $H'(u)= (1-u)^{-1}$ (4.14) and therefore
$$
u_k  = {1\over 1+\chi} = {k\over k+x_k}.
\eqno (4.15)
$$ 
 Then  
 $$
 T_k(x_k) = {\sqrt x_k\over \sqrt{2\pi k(x_k+{k})}}\, \left(1 + {x_k\over k}\right)^{k}
  \left(1 + {k\over x_k } \right)^{x_k} \,  (1+o(1)),\quad (k, x_k)_\chi\to\infty.
  \eqno (4.16)
  $$
  
  From the other hand, 
 using  identity \cite{Rio}
 $$
 T_k(x) = \sum_{(l_1,\dots,l_k)^*} {x^{l_1 + \dots +l_k}\over l_1! 
 \cdots l_k! \ 1^{l_1}\cdots k^{l_k} } = { x(x+1)\cdots (x+k-1)\over k!},
 \eqno(4.17)
 $$
 we see that (4.16) can be obtained  from (4.17) by  simple use of the Stirling formula (2.33) in the case when  $x_k= \chi k$, $k\to\infty$.
  It follows from relation (4.17) that 
$$
M_k(x) = {x\over x+k} \cdot { (x+k)!\over x!}.
$$
Regarding the asymptotic regime $1\ll x_k\ll k$, we can use the Stirling formula (2.33)
and write that 
$$
M_k(x)= \left({k+x\over e}\right)^k \left( {x+k\over x}\right)^{x-1/2}
\, \big(1 + o(1)\big).
$$
Remembering (4.15),  we get the following relation, 
$$
P(Z^{(x,u)}= k) = {M_k(x)\over k!} u^k (1-u)^x
$$
$$
= {k^k\over e^k k!} \left( {x+k\over x} \right)^{x-1/2} \left( {x\over k+x}\right)^x\big(1+o(1)\big) = {\sqrt x\over \sqrt{2\pi k(k+x)} } \big(1+o(1)\big), \quad k\to \infty.
\eqno (4.18)
$$
From another hand, the definition of $\s_Z^2$  means that 
$$
\s_Z = {\sqrt x u\over 1-u}\big(1+o(1)\big) = {k\over \sqrt x}\big(1+o(1)\big), \quad 1\ll x\ll k .
$$
Comparing this expression with the right-hand side of (4.18), we conclude that (cf. (2.21))
$$
P(Z^{(x,u)}=k) = {1\over \sqrt{2\pi } \s_Z}\big(1+o(1)\big),\quad (k,x_k)_0\to\infty
$$
and thus  that the Local Limit Theorem  holds for random variables $Z^{(x,u)}$.
This means that the restriction $m>8$ imposed 
on the Gamma distribution in subsection 3.1, or more generally, conditions (16), (1.7) and (1.8)  could be of  rather technical  character. 

\subsection{Concentration property of normalized sums}

Relations (1.13) is closely related with a law of large numbers for the random variable $Y_\l/\l$ as $\l\to\infty$ that is a known elementary fact. 
However, Theorem 1.1 gives more information about the limiting behavior of this variable. Indeed, 
given $y>0$, we deduce from the first asymptotic equality of (1.13)  that
$$
P\left( {1\over \l_k}  Y_{\l_k} >y\right) \le \left( {V_1\over y} \big(1+o(1)\big)\right)^k, \quad (x,k)_0\to \infty.
\eqno (4.19)
$$
The series of these probabilities converges for any $y>V_1$  and therefore by  the Borel-Cantelli lemma,
$$
P\left( \limsup_{k\to\infty} {Y_{\l_k}\over {\l_k}} \le V_1\right)= 1, \quad \l_k \gg k. 
\eqno (4.20)
$$

Regarding the moments of centered random variables  $\bar Y_\l= Y_\l - \l V_1$,
$$
\bar M_k(\l) = \bE \bar Y_\l^k
$$
 we can prove analog of Theorem 1.1. Indeed, one  can introduce auxiliary random variables $\bar Z^{(x,u)}$ by relation of the form (1.2),
 where $G(x,u)$ is replaced by 
 $$
 \bar G(x,u) = \sum_{j=0}^\infty \bar M_j(x) {u^j\over j!} = \exp\left\{ x \big( \bar H(u) - 1\big)\right\},
 $$
 where 
 $ \bar H(u) = H(u) - V_1u$. Then all computations of the proof of Theorem 1.1 can be literally repeated. 
 As a consequence, we can write in complete analogy with the second relation of (1.13) that
 for any $\vep>0$, 
 $$
 P\left( {1\over \l_k} |\bar Y_{\l_k} |>\vep\right)\le {1\over \vep^{2k} } \bar M_{2k}(\l_k) =  \left( {2V_2k\over e\l_k \vep^2}\big(1+o(1)\big)\right)^k,
 \quad 1\ll k \ll \l_k.
 \eqno (4.21) 
$$
This upper estimate implies convergence of $ Y_{\l_k}/\l_k$ to $V_1$  with probability 1 as $k\to\infty$ provided $k=o(\l_k)$ in this limit.

In the asymptotic regime when $\l_k = \chi k$, $k\to\infty$, we deduce from (1.15) the following version of (4.21),
$$
P\left( {1\over \l_k} |\bar Y_{\l_k} |>y\right)\le
 \left( { e^{\bar \Psi(\chi)} \over y} \big(1+o(1)\big)\right)^{2k},
 \eqno (4.22) 
$$
where 
$$
\bar \Psi(\chi) = {\bar H(u)- 1\over u \bar H'(u)} - 1 - \ln \bar H'(u)
$$
and 
$u$ is determined by equation 
$$
u \bar H'(u) = {1\over \chi}\, , \quad u>0.
$$

Upper bounds (4.20), (4.21) and (4.22) can be applied to a maximum of $n$ independent random variables 
$$
{1\over \l} Y_\l^{(i)}, \quad 
i=1, \dots , n
\eqno (4.23) 
$$
and their centered versions. Then one can obtain a number  of statements in the spirit of the Erd\H os-R\'enyi limit theorem \cite{ER}.
In particular, relation (4.22) will lead to the following upper bound
$$
P\left(
\max_{i=1, \dots, n}  
{1\over \l_k} \left|\bar Y^{(i)}_{\l_k} \right|>y\right)\le
n \left( { e^{\bar \Psi(\chi)} \over y} \big(1+o(1)\big)\right)^{2k}, \quad k\to \infty.
 \eqno (4.24) 
$$
It says that for any $C>0$, $k= C \log n$ and $\l_k= \chi k$, there exists $y= y(C,\chi)$ such that 
the superior  limit 
of random variables 
$$ 
T^{(n, \l_n)} = \max_{i=1, \dots, n} \,  
{1\over \l_n}  Y^{(i)}_{\l_n} 
\eqno (4.25) 
$$
remains bounded with probability 1 when $\l_n= \chi C  \log n$ as $n\to\infty$. 

The Erd\H os-R\'enyi limit theorem says that the maximum of $n$  random variables
$$
U^{(n,p)} =\max_{i=1, \dots, n} \,  {1\over p} \sum_{j=1}^p X_j^{(i)},
$$
where $\{X_j^{(i)}\}$ is a family of i.i.d. random variables converges, as $n\to\infty$  and $p= \tau \log n$
 to a non-random limit
$\alpha$  determined by $\tau $ and the exponential generating function of the moments of $X_j^{(i)}$. 
A generalization of this statement to the case of random variables $T^{(n,\l_n)}$ has been proved in \cite{K-02,K-03}.

\subsection{Maximal vertex degree of weighted random graphs}

Let us consider a family of i.i.d. random variables
$$
{\cal A}^{(n,\rho)}= \{ a_{ij}^{(n,\rho)}, \quad 1\le i< j\le n\},
$$
where $a_{ij}^{(n,\rho)}$ 
take values 1 and 0 with probability $\rho/n$ and $1 - \rho/n$. 
Real symmetric  random matrix  whose elements above the diagonal are given by (4.26) and by zero otherwise can be regarded as an adjacency matrix of a random graph with  $n$ vertices known as the Erd\H os-R\'enyi ensemble of random graphs \cite{Boll}. The weighted version is 
 given  by random symmetric matrices with the elements above the diagonal
$$
\left( A^{(n,\rho)}\right)_{ij} = a^{(n,\rho)}_{ij} \, X^{(i)}_j, \quad 1\le i< j\le n
$$
and zero on the diagonal. 
Then random variables
$$
D^{(i)}_\rho(n) = \sum_{j=1, \dots, n, \, j\neq i} \ a^{(n,\rho)}_{ij}\, X^{(i)}_j
\eqno (4.26)
$$
will play the role of the vertex degree of weighted random graphs. 
Regarding the random weights given by  i.i.d. $X^{(i)}_j$, $D^{(i)}_\rho(n)$ can be  regarded as a pre-limiting realization of random variables $Y^{(i)}_\l$ of (4.23). The only difference is that 
random variables $D^{(i)}_\rho(n)$ and $D^{(i')}_\rho(n)$ are not independent. However, this does not avoid the upper estimate of the 
deviation probabilities of the maximal vertex degree
$$
D^{(\max)}_\rho(n)= \max_{i=1, \dots, n} D^{(i)}_\rho(n).
$$
In analogy with (4.24), we can prove that 
given a sequence  $\rho_n= \kappa \log  n$, the following relation
$$
\lim_{n\to\infty} 
P\left( \left| D^{(\max)}_{\rho_n}(n)/\rho_n - V_1\right|
>s\right)= 0
\eqno(4.27)
$$
for any $s$ 
such that 
$$
s> \bar H'(u)\exp\left\{ {\bar H(u)-1 \over u \bar H'(u)}-{1\over 2}
\right\},
$$
where $u$ is determined by equation
$$
u\bar H'(u)={1\over \kappa}
$$
and 
$$
\bar H(u)= H(u)-u{\mathbb E}X_j^{(i)},  
\quad H(u)= \sum_{k=0}^\infty {u^k\over k!} \bE (X_j^{(i)})^k.
$$

If $\rho_n= \kappa_n \log n$ is such that $\kappa_n\to \infty$ as $n\to\infty$, then 
$$
\lim_{n\to\infty} D^{(\max)}_{\rho_n} (n)/\rho_n= V_1\quad {\hbox{with probability 1}}.
\eqno (4.28)
$$
Moreover, one can show that if a is a decreasing sequence $(s_n)_{n\in {\bN}}$ such that 
$$
s_n > \sqrt{ eV_2/\kappa_n} + V_1/n,
$$ 
then 
$$
P\left( \limsup_{n\to\infty}  \left| D^{(\max)}_{\rho_n}(n)/\rho_n - V_1\right|
>s_n\right)= 0.
\eqno (4.29)
$$

Relations (4.27), (4.28) and (42.9)  can be deduced from Theorem 1.1 with the help of standard 
 arguments of probability theory. We do not present the  proofs here. 

\subsection{Random matrices and even walks with multiple edges}

In this subsection we describe  one more situation when the moments of compound Poisson distribution can be useful in the random matrix theory.
Random real symmetric matrix of the Wigner ensemble \cite{W} is given by relation
$$
\left( W^{(n)}\right)_{ij} = {1\over \sqrt n} w_{ij}, 1 \le i < j \le n, \quad \left( W^{(n)}\right)_{ii} =0, 
$$
where ${\cal W}= \{ w_{ij}, 1\le i\le j\}$ is a family of joint independent identically distributed random variables. 
The eigenvalue distribution of the matrices $W^{(n)}$ can be studied with the help of the moments 
$$
\mu_l^{(n)} = \bE \left( {1\over n} \Tr \left(W^{(n)}\right)^l\right)= {1\over n} \sum_{i_1=1}^n\, 
\sum_{i_2, i_3, \dots, i_{k-1}} \bE\left( W^{(n)}_{i_1i_2} \cdot W^{(n)}_{i_2i_3}\cdots W^{(n)}_{i_{l}i_1}\right).
\eqno (4.30)
$$
The last sum of (4.30) can be regarded as a sum over all possible sequences $\CI_{l}^{(i_1)}= (i_2,i_3, \dots , i_{l})$ with corresponding weights. 
Adding the starting point $i_1$, this sequence can be represented as a multi-graph with 
 the set of vertices 
$\CV_n= \{1, 2, \dots, n\}$ and $k$ oriented edges. 

If one assumes that the probability distribution of $W_{ij}$ is symmetric, then  non-zero contributions to (4.30) 
are  given by sequences $\CI_{l-2}^{(i_1)}$, whose graphs have  vertices connected by an even number of edges. 
Thus one has to take  $l=2k$. 
A particular case of such an even sequence is given by 
$$
\tilde \CI_{2k}^{(i_1)}=(i_2, i_1, i_4, i_1,i_6,i_1,  \dots, i_1,i_{k}). 
$$ 
In the corresponding graph, we start with the vertex $i_1$, go to the vertex $i_2$, then return to $i_1$, then go to $i_4$, etcetera. 
In this construction, all variables $i_j$ with $j\neq 1$ differ from $i_1$.

It is not hard to see that the sum over all possible sequences $ \tilde \CI_{2k-2}^{(i_1)}$ with corresponding weights
$$
\CP(\tilde \CI_{2k}^{(i_1)}) = {1\over n^k} \bE \left( w_{i_1i_2} w_{i_2i_1} \, w_{i_1i_4} w_{i_4i_1}\,  \cdots\,  w_{i_1}w_{i_{2k}}w_{i_{2k}i_1}\right) 
$$
is given by expression of the form (1.2), also resembling (4.2A),
$$
\sum_{(i_2, i_4,i_6,  \dots,i_{k}), \, i_{2j} \neq i_1} \CP(\tilde \CI_{2k}^{(i_1)}) = {1\over n^k}\sum_{ (l_1, l_2, \dots, l_k)^*}^k
\prod_{i=1}^k {1\over l_i!} \left({ (n-1) W_{2i}\over i!}\right)^{l_i}   = {1\over n^k} M_k^{(W)}(n-1),
\eqno (4.31) 
$$
where we denoted $W_{2l}= \bE(W_{ij})^{2l}$, $l\in \bN$. 

Assuming that the function $H_W(u) = \sum_{j=0}^\infty u^{2j}  W_{2j}/(2j)!$ exists, we can use the results of Theorem 1.1
and say that 
$$
{1\over n^k} M_k^{(W)}(n-1) = \big( W_2(1+o(1))\big)^k, \quad 1 \ll k \ll n
\eqno (4.32)
$$
and that 
$$
{1\over n^k} M_k^{(W)}(n-1) = \left( e^{\Psi_W(\chi)}\big(1+o(1)\big)\right)^k, \quad k= \chi n, \ n\to\infty,
\eqno (4.33)
$$ 
where 
$
\Psi_W(\chi)$ and $u $ are  determined by  relations (1.16) and (1.17) with $H(u)$ replaced by $H_W(u)$. 
One can say that relation (4.32) considered with $k= Cn^{2/3}$ represents a kind of a   proof of a  known  result 
(see \cite {SS}, relation (4.29) and also \cite{K-01}), while the  estimate from below 
$
\mu_{2k}^{(n)} \ge  \left( e^{\Psi_W(\chi)}\big(1+o(1)\big)\right)^k
$ that follows from  (4.33) can be regarded as a new result. 
 


\end{document}